\RequirePackage{ifpdf}
\ifpdf 
\documentclass[pdftex]{sigma}
\else
\documentclass{sigma}
\fi

\begin{document}

\renewcommand{\PaperNumber}{103}

\FirstPageHeading

\renewcommand{\thefootnote}{$\star$}

\ShortArticleName{Geometric Linearization of Ordinary
Dif\/ferential Equations}

\ArticleName{Geometric Linearization \\ of Ordinary Dif\/ferential
Equations\footnote{This paper is a contribution to the Proceedings
of the Seventh International Conference ``Symmetry in Nonlinear
Mathematical Physics'' (June 24--30, 2007, Kyiv, Ukraine). The
full collection is available at
\href{http://www.emis.de/journals/SIGMA/symmetry2007.html}{http://www.emis.de/journals/SIGMA/symmetry2007.html}}}

\Author{Asghar QADIR}

\AuthorNameForHeading{A.~Qadir}

\Address{Centre for Advanced Mathematics and Physics,
National University of Sciences and Technology,\\
Campus of College of Electrical and Mechanical Engineering,\\
Peshawar Road, Rawalpindi, Pakistan}
\Email{\href{mailto:aqadirmath@yahoo.com}{aqadirmath@yahoo.com}}

\ArticleDates{Received August 13, 2007, in f\/inal form October
19, 2007; Published online November 06, 2007}

\Abstract{The linearizability of dif\/ferential equations was
f\/irst considered by Lie for scalar second order semi-linear
ordinary dif\/ferential equations. Since then there has been
consi\-derable work done on the algebraic classif\/ication of
linearizable equations and even on systems of equations. However,
little has been done in the way of providing explicit criteria to
determine their linearizability. Using the connection between
isometries and symmetries of the system of geodesic equations
criteria were established for second order quadratically and
cubically semi-linear equations and for systems of equations. The
connection was proved for maximally symmetric spaces and a
conjecture was put forward for other cases. Here the criteria are
brief\/ly reviewed and the conjecture is proved.}


\Classification{34A34; 34A26}

\renewcommand{\thefootnote}{\arabic{footnote}}
\setcounter{footnote}{0}

\section{Introduction}

Whereas there are standard methods to solve linear ordinary
dif\/ferential equations (ODEs), there are hardly any for {\it
nonlinear} ODEs. One could, of course, approximate the nonlinear
ODE by a~linear one (as in perturbation methods), but the key
features of the nonlinear equation that may be vital for the
phenomenon being modeled could easily be lost. As such, one needs
methods to obtain exact solutions of the equation. Further, if one
tries to improve the approximation by an iterative series the
convergence of the series would need to be proved. The same
problem is faced by numerical schemes for solution. Proof of
existence and convergence is provided by functional analytic
methods when they can be applied. On the other hand, some methods
{\it were} developed for specif\/ic classes of f\/irst order ODEs
and some second order ODEs were converted to linear form by
transformation of variables. Using symmetry considerations, Lie
had found that all f\/irst order ODEs could, in principle, be
converted to linear form by some appropriate change of independent
and dependent variables. (The inf\/initesimal generators of
symmetry form a Lie algebra \cite{Olver}.) He obtained criteria
for a second order semi-linear ODE to be such that it could be
converted to linear form by point transformations
\cite{Lie1880,Lie1883}. He did not go further to deal with systems
of ODEs or third order ODEs. His procedure uses transformations of
the dependent and independent variables of the ODE (or system of
ODEs) There has been work done on obtaining the symmetry classes
of systems of second order ODEs \cite{WafoSoh} and for the third
order ODE \cite{Leach}. Some explicit criteria for the
linearizability of single third order ODEs have also been provided
in the literature
\cite{Chern1937,Chern1940,Grebot1996,Grebot1997,Ibragimov,Meleshko}.

It is natural to ask what the connection is between symmetries in
geometry and for dif\/ferential equations. After all,
dif\/ferential equations ``live" on manifolds. More precisely, we
would like to relate the symmetry algebra of the ODEs to that of
the manifold. There had been some suggestions but there does not
appear to have been any attempt to formulate the connection
precisely till the turn of the century. The f\/irst such attempt
looked for the connection through the system of geodesic
equations\cite{Aminova}\footnote{There appears to be a Russian
preprint of 1991, but that is even less accessible than this
paper.}. It was used for the projective geometry of geodesics. (A
recent paper \cite{Bryant}, and references therein, detail
relations between projective algebra (point symmetries of a
particular class of second order ODEs) and algebra of isometries,
providing a~more complete background of geometry for the purpose.)
Independently, the same idea was developed and carried further to
construct the Lie algebra of the system of geodesic equations on
maximally symmetric spaces \cite{Feroze}. There a conjecture was
stated for all spaces. In this paper that conjecture will be
proved.

The plan of the paper is as follows. A brief review of the
geometrical  notation used will be presented in the next section.
The application of the geometric method for linearization of
second order ODEs will be presented in the subsequent section and
its extension to higher orders in Section~4. The theorem will be
proved in the subsequent section and a brief summary and
discussion presented in Section~6.

\section{Geometric notation used}

To explain the linearization procedure it is necessary to
introduce the geometric notation used. This will be presented
f\/irst.

The following are well-known and can be found in text books. We
use the Einstein summation convention that repeated indices are
summed over the entire range of the index. Thus, $A^aB_a$~stands
for $\sum_{a=1}^n A^aB_a$. The metric tensor will be represented
by the symmetric (non-singular) matrix $g_{ij}$ and its inverse by
$g^{ij}$. The Christof\/fel symbols are given by
\begin{gather*}
\Gamma^i_{jk}=\frac12 g^{im}(g_{jm,k}+g_{km,j}-g_{jk,m}),
\end{gather*}
where $,k$ stands for partial derivative relative to $x^k$, etc.
The Christof\/fel symbols are symmetric in the lower pair of
indices, $\Gamma^i_{jk}=\Gamma^i_{kj}$.

In this notation, the system of $n$ geodesic equations is
\begin{gather*}
\ddot x^i+\Gamma^i_{jk}\dot x^j\dot x^k=0,\qquad i,j,k=1,\ldots,n,
\end{gather*}
where $\dot x^i$ is the derivative relative to the arc length
parameter $s$ def\/ined by $ds^2=g_{ij}dx^i dx^j$.

The Riemann tensor is def\/ined by
\begin{gather*}
R^{i}_{jkl}=\Gamma^i_{jl,k}-\Gamma^i_{jk,l}+\Gamma^i_{mk}\Gamma^m_{jl}-\Gamma^i_{ml}
\Gamma^m_{jk}, 
\end{gather*}
and has the properties $R^i_{jkl}=-R^i_{jlk}$,
$R^i_{jkl}+R^i_{klj}+R^i_{ljk}=0$ and
$R^i_{jkl;m}+R^i_{jlm;k}+R^i_{jmk;l}=0.$ The Riemann tensor in
fully covariant form $ R_{ijkl}=g_{im}R^m_{jkl}$ satisf\/ies
$R_{ijkl}=-R_{jikl}.$

\section{Linearizability of second order ODEs}

A general quadratically semi-linear system of ODEs is of the form
\begin{gather}
{x^a}^{\prime \prime}+ \alpha^a + \beta^a_b{x^b}^{\prime} +
\gamma^a_{bc}{x^b}^{\prime}{x^c}^{\prime} =0. \label{(4)}
\end{gather}
The system of ODEs will be said to be of {\it geodesic form} if
$\alpha^a=\beta^a_b=0$. It can be regarded as a~system of geodesic
equations if there exists some metric tensor for which the
Christof\/fel symbols satisfy $\Gamma^a_{bc}=\gamma^a_{bc}$. In
this case we can construct the curvature tensor corresponding to
the given coef\/f\/icients. It is easier to see the procedure
adopted for the case of two dependent variables. Call the
variables~$x$,~$y$. Then we have 6 arbitrary functions of these
variables appearing as coef\/f\/icients in the equations
\begin{gather*}
x''=a(x,y)x'^2+2b(x,y)x'y'+c(x,y)y'^2,\nonumber\\
y''=d(x,y)x'^2+2e(x,y)x'y'+f(x,y)y'^2,
\end{gather*}
the Christof\/fel symbols in terms of these coef\/f\/icients are
\begin{gather*}
\Gamma_{11}^1=-a, \qquad \Gamma^1_{12}=-b, \qquad
\Gamma_{22}^1=-c,\qquad
\Gamma_{11}^2=-d, \qquad \Gamma^2_{12}=-e, \qquad \Gamma_{22}^2=-f.
\end{gather*}
The linearizability criteria are \cite{Qadir}
\begin{gather}
a_y-b_x + be -cd = 0 ,\qquad  b_y -c_x+ (ac - b^2) + (bf - ce) = 0 ,\nonumber\\
d_y-e_x - (ae - bd) - (df - e^2) = 0 , \qquad (b + f)_x = (a +
e)_y , \label{(7)}
\end{gather}
with constraints on the metric coef\/f\/icients: $g_{11}=p$,
$g_{12}=g_{21}=q$, $g_{22}=r$,
\begin{gather*}
p_x = -2(ap + dq) , \qquad q_x = -bp - (a + e)q - dr , \qquad r_x = -2(bq + er) ,\nonumber\\
p_y = -2(bp + eq) , \qquad q_y = -cp - (b + f)q - er , \qquad r_y
= -2(cq + fr) .
\end{gather*}
The compatibility of this set of six equations gives the above
four linearization conditions \eqref{(7)}. {\it The remarkable
fact is that the linearization criteria are simply} $R^i_{jkl}=0$!
For a system of three equations we get eighteen such equations. In
general there are $n^2(n+1)/2$, $n\ge 2$. One now obtains the
required linearizing transformation by regarding the variables in
which the equation is linearized as Cartesian, thus having
$g_{11}=g_{22}=1$ and $g_{12}=g_{21}=0$, and looking for the
coordinate transformation yielding the original metric
coef\/f\/icients \cite{Qadir}. This can be done by using complex
variables. The details can be found in the reference, where
examples are also given.

Following Aminova and Aminov \cite{Aminova} and projecting the
system  $n$-dimensional system down to a system of $(n-1)$
dimensions, the geodesic equations become \cite{Mahomed}
\begin{gather}
{x^a}^{\prime\prime}+A_{bc}{x^a}^{\prime}{x^b}^{\prime}{x^c}^{\prime}+B^a_{bc}{x^b}^{\prime}
{x^c}^{\prime}+C^a_b{x^b}^{\prime}+D^a=0,\qquad
a=2,\ldots,n,\label{(9)}
\end{gather}
where the prime now denotes dif\/ferentiation with respect to the
parameter $x^1$ and the coef\/f\/icients in terms of the
$\Gamma^a_{bc}$s are
\begin{gather*}
A_{bc}=-\Gamma^1_{bc}, \qquad
B^a_{bc}=\Gamma^a_{bc}-2\delta^a_{(c}\Gamma^1_{b)1},
\nonumber\\
C^a_b=2\Gamma^a_{1b}-\delta^a_b\Gamma^1_{11}, \qquad
D^a=\Gamma^a_{11},\qquad
a,b,c=2,\ldots,n.
\end{gather*}
System \eqref{(9)} is cubically semi-linear. When $n=2$, writing
$a-2e=h$, $f-2b=g$, we get
\begin{gather}
y''+c(x,y)y'^3-g(x,y)y'^2+h(x,y)y'-d(x,y)=0. \label{(11)}
\end{gather}
This is the general form of linearizable second order semi-linear
ODE found by Lie, with the coef\/f\/icients re-named. To be
linearizable it must satisfy the linearizability
conditions{\samepage
\begin{gather*}
h_y+2e_y-b_x+be-cd=0, \qquad b_y-c_x+hc+ec+b^2-bg=0,\nonumber\\
-d_y+e_x+eh+e^2+bd-dg=0,\qquad  3b_x-3e_y-g_x-h_y=0, 
\end{gather*}
{\it which are simply the Lie criteria for a scalar ODE\/}!}

One can use the same procedure to go on to systems of two (or
more) cubically semi-linear ODEs by taking $n \ge 3$ in
\eqref{(4)}. This has been done~\cite{Mahomed} and examples of the
two dimensional case are also given there.

It is known that there is a unique class of linearizable scalar
second order semi-linear ODEs, namely that which satisf\/ies the
Lie criteria. One may  wonder if the extension provided for
systems will give a unique class as well. The answer is ``No'', as
it is known that for $2$-dimensional systems there are f\/ive
classes \cite{WafoSoh} that are linearizable and only one has been
provided by the above procedure. How is it that the case of the
single equation is fully covered by the above procedure but the
systems are not fully covered?

Note that the projection procedure applied to an arbitrary system
of $n$ quadratically semi-linear second order ODEs reduces it to
$n-1$ cubically semi-linear second order ODEs. (Of course there is
an enormous computational complication arising. To construct the
metric tensor from the Christof\/fel coef\/f\/icients an algebraic
computational code has been written \cite{Fredericks}.) Now
observe that in projecting down from the system of $n$ variables
to $n-1$ variables, the Christof\/fel symbols are reduced from
$n^2(n+1)/2$ by $n$, to give $(n-1)n(n+2)/2$ independent
coef\/f\/icients. Since we now have $n-1$ equations, each with its
own cubic function, there are $4(n-1)$ coef\/f\/icients for the
reduced system. If the number of coef\/f\/icients left after
losing $n$ equals the number of coef\/f\/icients of the reduced
system, we can determine one set of coef\/f\/icients in terms of
the other. The two expressions are obviously equal for $n=2$ and
the former is greater than the latter for $n>2$. As such, the
coef\/f\/icients of the cubic system can be determined uniquely in
terms of the quadratic system for $n=2$, i.e.\ for a scalar
cubically semi-linear system. For larger systems there will be
inf\/initely many ways to write the former in terms of the latter.
This is why there is a unique solution to the linearizability
problem for the scalar equation and many solutions for the
systems! It can be seen that \eqref{(9)} is {\it not} the most
general form for the system as one should have
$A^a_{bcd}{x^b}^{\prime}{x^c}^{\prime}{x^d}^{\prime}$ in place of
the second term for full generality. For $n=2$ there can be no new
term coming from this more general coef\/f\/icient, but for $n \ge
3$ new terms will come in. {\it This is where the extra classes
that are missed for $n \ge 3$ will come from.}

\section{Higher order systems}

A procedure to extend the scalar second order case to the third
order is to dif\/ferentiate the general second order ODE
\eqref{(11)}. One can now replace the second order derivative term
in the equation by using the original second order ODE. Thus
\cite{MahomedQadir} one either gets a third order ODE which is
quintic in the f\/irst derivative
\begin{gather*}
y'''-\alpha y'^5+\beta y'^4-\gamma y'^3+\delta y'^2-\epsilon
y'+\phi=0, 
\end{gather*}
or a third order ODE with a term in the second order that has a
coef\/f\/icient quadratic in the f\/irst derivative
\begin{gather*}
y'''-(A(x,y)y'^2-B(x,y)y'+C(x,y))y''+D(x,y)y'^4+E(x,y)y'^3\nonumber\\
\qquad{}-F(x,y)y'^2+G(x,y)y'-H(x,y)=0.
\end{gather*}
It turns out that the same result is obtained in either case.
These equations are not in the classes considered in the
literature earlier. It must be added that the class of equations
obtained by this procedure also {\it does not} contain the classes
earlier discussed in the literature . Instead of following this
procedure, one could have gone in the reverse direction. First
dif\/ferentiate the system and then project it. Surprisingly, it
turns out that the result is again the same in either case. Some
examples of linearizable third order equations are given in
\cite{MahomedQadir}, one of which overlaps with the class of
Meleshko \cite{Meleshko}, but none of which overlaps with that of
Ibragimov and Meleshko~\cite{Ibragimov}.

It should be noted that if one tries to reduce from the quintic,
it is not possible to reduce to the quartic. Further, the cubic
cannot be resolved directly as there are problems of
non-uniqueness of the solutions of the equations for the
coef\/f\/icients required to linearize the third order cubic ODE.
It would be useful to f\/ind a canonical way to solve the
equations, or to prove that all solutions are equivalent. One
needs to understand how a complete classif\/ication of the various
classes can be obtained and to provide explicit criteria and
constructive procedures for the solutions. In particular, there
are only three classes possible for the third order ODEs
linearizable by point transformations \cite{Leach},
$y^{\prime\prime\prime}(x)=0$, $y^{\prime\prime\prime}(x)=y(x)$
and $y^{\prime\prime\prime}(x)=\alpha (x)y(x)$. Two were provided
by Chern \cite{Chern1937,Chern1940} and the third by Neut and
Petitot \cite{Neut}, and later by Ibragimov and Meleshko
\cite{Ibragimov}. How can there then be {\it another} class? The
point is that the method just outlined uses a projection and
dif\/ferentiation apart from the point transformations. It
guarantees only two arbitrary constants and not three while the
other procedure guarantees three.

Can one go further? On dif\/ferentiating the third order ODE (in
the quintically semi-linear form) and using the second order ODE,
just as in the case of the third order cubically semi-linear ODE,
there is no guarantee that the coef\/f\/icients of the fourth
order ODE would be resolvable in terms of the coef\/f\/icients of
the original second order ODE. It turns out that it is
\cite{QadirMahomed} resolvable. However, as yet there is no proof
that the other forms of the equation, involving the third and
lower orders or the second and f\/irst orders, would yield the
same result. It certainly is not clear that the projection and
dif\/ferentiation procedures would commute. It is not even clear
how that would be checked. One would need a system of third order
equations to be able to check this result. Though the system of
third order equations that are cubically semi-linear in the the
f\/irst derivative has been computed \cite{Naeem}, it is highly
complicated and one cannot yet make any statement about the
commutation of the dif\/ferent procedures.

\section[Proof of the Feroze-Mahomed-Qadir conjecture]{Proof of the Feroze--Mahomed--Qadir conjecture}

In \cite{Feroze} maximal curvature spaces were considered and it
was found that for f\/lat spaces of dimension $n$ the symmetry
algebra is $sl(n+2,{\mathbb R})$, for positive curvature
$so(n+1)$, or for negative curvature a corresponding non-compact
version. There was no result for spaces that are not maximally
symmetric. However, one can imagine some invariant procedure to
split the space into f\/lat, positive or negative curvature
spaces, or those with no symmetry. Considering each of these
subspaces it seemed clear that the corresponding direct sums
should provide the full symmetry algebra. In case there is no
f\/lat subspace, there will still remain the geodesic
re-parametrization by translation and scaling. This is what gives
the extra $2$ in $sl(n+2,{\mathbb R})$. Where there is no part of
the special linear algebra the dilation algebra, $d_{2}$, will
replace it. On this basis a more general conjecture was stated.
That conjecture is now proved as a theorem.

\medskip

\noindent {\bf Theorem.} {\it For a space of non-zero curvature
with isometry algebra $h$  the symmetry algebra of the geodesic
equations is $h\oplus d_{2}$  provided that there is no section of
zero curvature. If there is an $m$-dimensional maximal
cross-section of zero curvature, $M$, and the symmetry algebra of
the orthogonal subspace, $M^{\perp }$  is $h_{1}$,  the symmetry
algebra of the geodesic equations will be $h_{1}\oplus
sl\left(m+2,{\mathbb R} \right)$. }

\begin{proof} First consider a base manifold $N$, which has a maximal
f\/lat section $M$ of dimension $m\ne 0$ such that $N=M\times
M^{\perp}$, where the orthogonal subspace, $M^{\perp}$, has no
f\/lat section. Consider a metric tensor ${\bf g}$ on $N$ and let
the isometry algebra on $M^{\perp}$ be~$h$. Thus the isometry
algebra of $N$ will be $so(m) \oplus_s \oplus h$, where $\oplus_s$
denotes the semi-direct sum, i.e.\ a ``non-commutative direct
sum''.

Construct a f\/ibre bundle on $N$ with the f\/ibre given by
${\mathbb R}$ and take the position vector ${\bf x}$ to depend on
the arc length, or geodetic, parameter along the f\/ibre $s$,
i.e.\ such that ${\bf g.\dot x \dot x}$ is unity. If one def\/ines
the Lagrangian, as a function of $\bf x$ and ${\bf \dot x}$ to be
this quantity, it will possess at least one additional symmetry
generator, $\partial /\partial s$, apart from the isometry
algebra. The Lagrangian is not, in general, invariant under the
generator $s \partial / \partial s$ (unless it has the value
zero). However, the Euler--Lagrange equations for this Lagrangian
{\it will} be invariant under this additional generator.

Consider the geodesic equations for $M$ in the f\/ibre bundle over
it. They will be the vector free-particle equations and hence will
have the symmetry algebra $sl\left(m+2,{\mathbb R} \right)$. The
remaining part will retain its isometry algebra. Thus, in this
case, the full symmetry algebra is $h_{1}\oplus
sl\left(m+2,{\mathbb R} \right)$.

Now consider the case that there is no f\/lat section of the
manifold, $M=\varnothing$ so $m=0$. In this case the isometry
algebra will simply be~$h$. There will remain the two symmetry
generators for re-parametrization of the geodetic parameter by
translation and re-scaling, $\partial /\partial s$ and $s\partial
/ \partial s$. They form a two-dimensional dilation
algebra,~$d_2$.
\end{proof}

This theorem provides the basis for using geometric procedures to
obtain linearizable forms of (systems of) ODEs, through the link
between the symmetries of the manifold and the system of geodesic
equations on it. In particular, one requires a f\/lat space for
the system of geodesics to be linearizable. The geometry of the
projective space need not be f\/lat but it must be imbeddable in a
f\/lat space. While it was used as a conjecture, one really needed
a rigorous proof of the result.

\section{Conclusion}

In this paper a review has been provided of the geometric approach
to symmetry analysis using the connection between symmetries of
geometry and of dif\/ferential equations provided by the system of
geodesic equations, noted by Aminova and Aminov~\cite{Aminova} and
separately by Feroze et al.~\cite{Feroze}. This leads to a new
approach to investigate the linearizability of ODEs. A bonus of
this approach is that the solution of the linearized equations is
obtainable by constructing the transformation of coordinates from
Cartesian to the metric tensor given by the geodesic
equations~\cite{Qadir}. Using this approach and the projection
procedure of Aminova and Aminov~\cite{Aminova} one can re-derive
the Lie criteria for a scalar equation and extend them to a system
of equations~\cite{Mahomed}. One can also extend to the third and
higher order scalar equations~\cite{MahomedQadir} and to the third
order system of equations~\cite{Naeem}. There is much work that
needs to be done in this direction. For one thing, the systems
have many classes of which one segment has been dealt with so far.
A~complete analysis and a method to f\/ind solutions of the
remaining classes would be very important. Again, the higher order
equations can be linearized by using the Lie procedure as done by
Ibragi\-mov and Meleshko \cite{Ibragimov} or by the procedure of
Neut and Petitot \cite{Neut}. To what extent is there an overlap
between the two procedures? For that matter, does this work
subsume the work of Grebot \cite{Grebot1996,Grebot1997} or not?
The relation between the various methods and approaches being used
is needed for a proper understanding and classif\/ication of the
procedures.

A conjecture that was at the base of the connection between
geometry and symmetry analysis found by Feroze et al \cite{Feroze}
was stated and proved as a theorem. The proof given demonstrates
that there is at least one extra symmetry generator of the
Lagrangian apart from the isometry algebra, and that the symmetry
algebra of the geodesic equations contains the symmetry algebra of
the Lagrangian. The Lagrangian can have additional symmetries
beyond the one extra symmetry generator of geodetic
re-parametrization by translation due to ``mixing'' between the
additional dimension of the f\/ibre bundle and the other f\/lat
dimensions. It would be of interest to use the above approach in
proving the conjecture to determine the symmetry algebra of the
Lagrangian. The problem is that whereas for the Euler--Lagrange
equations the ``mixing of the two symmetries of the f\/ibre with
the f\/lat base space'' is apparent and is maximal, that is not
true for the Lagrangian. For one thing, there is generally only
one symmetry for the f\/ibre and the mixing along a section of the
f\/ibre is {\it not} maximal and is not apparent. How, and to what
extent, the extra symmetries from the f\/ibre come into the
Lagrangian needs to be determined rigorously.

\subsection*{Acknowledgements}

I am grateful to NUST for travel support and to the organizers of
the Symmetry-2007 and the International Mathematical Union for
local support and hospitality at the conference where this paper
was presented. I am also grateful for useful comments to
Professors Leach, Mahomed, Meleshko and Popovych. Thanks also to
DECMA and CAM of the of Wits University, Johannesburg, South
Africa for support at the University where the paper was
completed.

\pdfbookmark[1]{References}{ref}
\LastPageEnding

\end{document}